\documentclass[11pt, leqno]{amsart}

\usepackage{graphicx}

\usepackage{amsfonts,amsmath,oldgerm,amssymb,amscd,mathrsfs,bbm,palatino}

\usepackage[all]{xy}

\newcommand{\ol}{\overline}

\newcommand{\wt}{\widetilde}

\newcommand{\rar}{\rightarrow}

\newtheorem{proposition}{Proposition}

\newcommand{\Z}{\mbox{$\mathbb Z$}}

\newcommand{\nind}{\noindent}

\newcommand{\Q}{\mbox{$\mathbb Q$}}

\newcommand{\C}{\mbox{$\mathbb C$}}       

\begin{document}

\title{Fundamental Group of some Genus-2 fibrations and applications}

\author{R.V. Gurjar}

\author{Sagar Kolte}

\noindent
\footnote{Mathematics Subject Classification 2010, {14F35, 14H30}}


\begin{abstract}

We will prove that given a genus-2 fibration $f: X \rightarrow C$ on a smooth projective surface $X$
such that $b_1(X)=b_1(C)+2$, the fundamental group of $X$ is almost isomorphic to
$\pi_1(C) \times \pi_1(E)$, where $E$ is an elliptic curve. We will also verify
the Shafarevich Conjecture on holomorphic convexity of the universal cover of surfaces $X$ with genus-2 
fibration $X\rar C$ such that $b_1(X)>b_1(C)$.

\end{abstract}

\maketitle

\section*{Introduction}

\nind
Let $X$ be a smooth projective surface and $C$ a smooth projective curve over
$\C$. Suppose that there exists a morphism $f:X \rar C$ such that a general
fiber is a smooth curve of genus 2. It was proved in \cite{GP1} that for a genus
two fibration with irreducible singular fibers and ample canonical divisor,
the image of the map $i: \pi_1(F) \rightarrow \pi_1(X)$ is finite, where $F$ is
a general fiber. This was later generalized to hyperelliptic fibrations in
\cite{GPP}. Using this the Shafarevich Conjecture for such fibrations was
proved. In this paper we will prove the following result.\\

\nind
{\bf Theorem.} {\it Let $f:X\rar C$ be a relatively minimal genus-$2$ fibration on a smooth
projective surface $X$ onto a smooth projective curve $C$. Suppose that $b_1(X)>b_1(C)$. Then we have the following assertions.\\

(1) If $f$ is not a $C^{\infty}$-fiber bundle then one of the following holds:\\

(i) $$\pi_1(X)\cong \pi_1(C)\times ({\bf Z}\oplus{\bf Z})$$ where ${\bf
Z}\oplus{\bf Z}$ is the fundamental group of an elliptic curve.\\

(ii) $\pi_1(R)\times \pi_1(E_0)$ maps homomorphically to $\pi_1(X)$ with image of index at most $2$, where $R$ is a smooth projective curve with a two-to-one morphism to $C$ and $E_0$ is a smooth elliptic curve. In this case, $\pi_1(X)\cong\Gamma_1\times\pi_1(E_0)$, where $\Gamma_1$ is the image of the map $\pi_1(R)\rar\pi_1(X)$. \\

(2) Suppose that $f$ is a $C^{\infty}$-fiber bundle.\\
If $C$ is a rational curve then $X$ is isomorphic to $C\times F$.\\
If $C$ is an elliptic curve then $X$ is isomorphic to a quotient of a trivial bundle $C'\times F$, where $C'$ is an elliptic curve and $F$ is a smooth curve of genus $2$.\\

(3) If $b_1(X)=b_1(C)+4$ then $X$ is isomorphic to $C\times F$. In fact, if the genus of $C$ is either $0$ or $>2$ then $f$ is itself a trivial bundle.}\\

{\it Remarks.} (a) If $f$ is a $C^{\infty}$-fiber bundle then we have the well-known exact sequence 
$$\cdots\rar\pi_2(X)\rar\pi_2(C)\rar\pi_1(F)\rar\pi_1(X)\rar\pi_1(C)\rar (1).$$
If further the genus of $C>0$ then $\pi_i(C)=(0)$ for $i>1$.\\

\nind
(b) Part (1) of the theorem is the main result of this paper. The first part of (3)
in the theorem is due to A. Beauville \cite{Beu}. Presumably, the second part of (3) is well-known to experts. We prove in the preliminaries (see, Proposition) that if $f$ is a relatively minimal genus-2 fibration which has at least one singular fiber then $b_1(X)\leq b_1(C)+2$.\\
If on the other hand, if  $b_1(X)=b_1(C)+4$ then $X$ is isomorphic to $C\times F$.\\

In the last section we will use the theorem to verify of the Shafarevich Conjecture for genus-2 fibrations
such that $b_1(X)>b_1(C)$.\\

\nind
{\bf Corollary 1.} {\it With the notation and assumptions of the theorem if $b_1(X)>b_1(C)$ then the  universal covering space $\widetilde X \rar X$ is holomorphically convex.}\\

\nind
{\bf Corollary 2.} {\it With the notation and assumptions of the theorem if $b_1(X)>b_1(C)$ then the
second homotopy group $\pi_2(X)$ is a free abelian group.}\\

\nind
{\it Acknowledgements.} The authors would like to thank the referee for reading the paper carefully, pointing out some inaccuracies and for  suggesting improvements in the presentation.\\

\section{Preliminaries.}

\nind
In this paper we will deal only with complex algebraic or complex analytic varieties.\\

\nind
We will start with some useful general results about fibrations on smooth projective surfaces.\\
By a component of a (possibly reducible and non-reduced) curve $B$ we mean an irreducible component of $B$.\\
By a rational (elliptic) curve we mean a possibly singular irreducible curve whose desingularization has genus $0$ (resp. $1$).\\
A smooth projective rational curve with self-intersection $-n$ on a smooth projective surface is called a $(-n)$-curve.\\
Recall that an automorphism $\sigma$ of finite order of a local analytic domain
$R$ is a {\it pseudo-reflection} if there is a minimal generating set
$x_1,x_2,\ldots,x_r$ for the maximal ideal of $R$ such that $\sigma(x_i)=x_i$
for $1\leq i\leq r-1$ and $\sigma(x_r)=\omega x_r$ for a suitable root of unity
$\omega$. We will implicitly use the well-known result of Shephard-Todd,
Chevalley that if $R$ is the power series ring and $G$ is a finite group of
automorphisms of $R$ generated by pseudo-reflections then the ring of invariants
$R^G$ is again isomorphic to a power series ring \cite{C}.\\

\nind
Let $\varphi:Y\rar B$ be a morphism on a smooth projective surface $Y$ onto a smooth projective curve $B$ of genus $g$ such that a general fiber $F$ of $\varphi$ is irreducible. Let $m_1F_1,\ldots,m_rF_r$ be all the singular fibers of $\varphi$ with multiplicities $m_i>1$ for each $i$. Then there is a short exact sequence
$$\pi_1(F)\rar\pi_1(Y)\rar\Gamma\rar (1)$$
where $\Gamma\cong <a_1,b_1,..,a_g,b_g,c_1,..,c_r|[a_1,b_1]...[a_g,b_g]c_1c_2...c_r=1,c_1^{m_1}=...=c_r^{m_r}=1>$ (\cite{X}, Lemma 2). Using this we see that the following natural sequence of homomorphisms is exact

$$H_1(F;{\bf R})\rar H_1(Y;{\bf R})\rar H_1(B;{\bf R})\rar (0).$$

Thus, $b_1(Y)\leq b_1(B)+b_1(F)$. On the other hand, A. Beauville has
proved in \cite{Beu} that if equality occurs then $Y$ is birational to $B\times
F$. In this case $\pi_1(Y)$ is isomorphic to $\pi_1(B)\times\pi_1(F)$.\\

We have the following useful result which proves part (3) of the theorem.

\begin{proposition}
 Let $f: W \rightarrow B$ be a relatively minimal genus-2 fibration on a smooth
projective surface $W$ onto a smooth curve $B$.\\

(a) If $f$ has a singular fiber then $b_1(W) \leq b_1(B) +2$.\\

(b) If $b_1(W)=b_1(C)+4$ then $X$ is isomorphic to $C\times F$.\\
If the genus of $C$ is either $0$ or $>2$ then $f$ is a trivial bundle.
\end{proposition}

\begin{proof}
(a) We will first prove the result when genus of $B > 2$ and reduce the
general case to this situation.\\
Let genus $B \textgreater 2$. Suppose $b_1(W)=b_1(B)+4$. By Beauville's result
$W$ is birational to a product $B \times F$ where $F$ is a smooth projective
surve of genus 2. There is a surface $\wt W$ obtained by a sequence of
blow ups from $W$ such that the induced birational map $\tau: \wt W
\rightarrow B \times F$ is a morphism. Since $W \rightarrow B \times F$ is
defined outside a finite set of points, a general fiber $F_0$ of $f$ is also a
general fiber of $\tilde f: \wt W \rightarrow B$, and the image $\tau(F_0)$
has to be a fiber of the projection $B \times F \rightarrow B$ since genus
$B\textgreater 2$. This shows that the birational map $W \rightarrow B \times F$
sends fibers of $f$ to fibers of the projection $B \times F \rightarrow B$.\\

Let $F_s$ be a singular fiber of $f$. Then both $F_s$ and its inverse image in
$\wt W$ are
unions of rational and elliptic components. But such a fiber in $\wt W$
cannot map onto a genus 2 curve $F$ which is a fiber of $B \times F \rightarrow
B$.\\
This contradiction shows that $ b_1(W) \leq b_1(B) +2$.\\

Now we consider the general case. We can find a finite ramified map $\wt B \rightarrow B$ which
is {\'e}tale at all points in $B$ over which a singular fiber of $f$ occurs and
such that genus $\wt B \textgreater 2$. Consider the fiber product $\wt W:= W
\times_B \wt B$. Then $\wt W$ is smooth and has an induced genus-2 fibration
$\wt f : \wt W \rightarrow \tilde B$. Clearly $\wt f$ is relatively
minimal and $\wt f$ has at least one singular fiber. By the previous case,
the image $H_1(F;\Q) \rightarrow H_1(\wt W;\Q)$ has rank $\leq 2$, where $F$ is
a
general fiber of $\wt f$. Then $H_1(F;\Q) \rightarrow H_1(W;\Q)$ also has
rank $\leq 2$ for a general fiber $F$ of $f$.\\

(b) Now assume that $b_1(W)=b_1(C)+4$. By (a) $f$ is a $C^{\infty}$-fiber bundle
and $X$ is birational to the product $C\times F$.\\
If $C$ is rational then by Teichm{\"u}ller theory $f$ is a trivial bundle.\\

Assume that genus $C>0$. We will first prove that the birational (rational) map
$W\rightarrow C\times F$ is a morphism.\\
If this is not true then let $\wt W\rightarrow X$ be obtained by a smallest
succession of blowing ups so that the induced birational map $\wt W\rightarrow
C\times F$ is a morphism. Let $E$ be the last exceptional curve obtained in the
blowing ups. Then the image of $E$ in $C\times F$ is a curve. This curve cannot
map onto the base $C$ since genus $C>0$. This is a contradiction. Hence we have
a birational morphism $W\rightarrow C\times F$. 
It is clear that $W$ cannot contain any rational curve since $f$ is a
$C^{\infty}$-bundle and $C$ is non-rational. Hence this birational morphism is
also finite. It follows that this morphism is an isomorphism.\\
Hence $W$ is isomorphic to $C\times F$.\\

Finally, assume that genus $C>2$. Then the proof of part (a) shows that the fibers of $f$ are mapped to fibers of the projection $C\times F\rightarrow C$. Hence $f$ is a trivial bundle.\\
This proves part (b) of the proposition.\\ 
\end{proof}

{\it Remark.} In (b) above it is not clear if $f$ is itself a trivial bundle in case $C$ is of genus $1$ or $2$.\\

In view of the proposition we can now assume for the rest of the proof that $$b_1(X)=b_1(C)+2$$

\nind
We will need to use precise knowledge of singular fibers of a genus-2 fibration
$X\rar B$ on a smooth projective surface $X$ from \cite{NU}.\\
For the purpose of calculating the fundamental group of $X$ in this paper the
singular fibers which contain an elliptic curve need a careful scrutiny.\\
Firstly, every singular fiber contains an irreducible component which occurs
with multiplicity $1$ in the scheme-theoretic fiber. There are only two types of
singular fibers which contain an elliptic curve occurring with multiplicity $2$
in the fiber, viz. those of numerical type (12) on pages 155 and 159 in
\cite{NU}. Next, suppose that a singular fiber contains an elliptic curve
(which is the only case of real interest in this paper by Lemma 1 below). In the
Step 1 in the proof we eliminate most types of singular fibers containing a reduced elliptic curve as a component. After this step
we have the situation that the union of rational curves in any singular fiber is
a disjoint union of is simply-connected curves. The union of rational curves can
be contracted to finitely many normal singular points on a normal compact
complex surface $T$. An easy application of Van Kampen's theorem shows that the
induced homomorphism $\pi_1(X)\rar\pi_1(T)$ is an isomorphism. This observation
will be rather useful in some of the arguments.\\

In our proof we will first deal with the cases when at least one singular fiber contains a reduced elliptic curve and then deduce the result in the two cases when only multiple elliptic curves occur in singular fibers.\\

\section{Proof of the Theorem.}

\nind
We now assume that $b_1(X)=b_1(C)+2$.\\

\nind
\nind
We have an induced exact sequence of abelian varieties and homomorphisms
$$Alb(F)\rar Alb(X)\rar Alb(C)\rar (0)$$ and the image of the homomorphism
$Alb(F)\rar Alb(X)$ is a $1$-dimensional abelian subvariety of $Alb(X)$ denoted
by $E'$. This will be essentially the elliptic curve $E$ in the statement of the
theorem.\\

\nind
By Poincar{\'e}'s complete reducibility we can find an abelian subvariety
$S\subset Alb(X)$ such that $E'\cap S$ is finite and $S$ maps onto $Alb(C)$ with
a finite kernel. Now $Alb(X)/S$ is an elliptic curve $E_0$. This gives a
surjective morphism $X\rar E_0$. Let $g:X\rar E$ be the Stein factorization of
this morphism.\\

{\bf Proof of Part (1) of the Theorem.}\\

\nind
Now we will assume that $f$ is not a $C^{\infty}$-fiber bundle. We will use the knowledge of possible singular fibers in a genus-2
fibration.\\

\nind
In \cite{NU} and \cite{O} the authors have given a complete classification of
the singular fibers that can occur in a genus-2 fibration. There are $44$
different possible types of singular fibers. Any component of each fiber is
either a rational curve or an elliptic curve. The next result says that many of
these fibers cannot occur in our situation.\\

\nind
{\bf Lemma 1.} {\it Every singular fiber of $f$ contains an elliptic curve.\\
If a singular fiber of $f$ is irreducible then it is an elliptic curve with exactly one ordinary node and no other 
singular point.}\\

\nind
{\bf Proof.} Since $f$ is flat, every scheme-theoretic fiber has arithmetic genus $2$.\\
Suppose that a singular fiber $F_0$ has all its components only rational curves. Then the image of $F_0$ in $Alb(X)$ is a point. By continuity, all the fibers of $f$ are mapped to points in $Alb(X)$. This implies that we have an isomorphism $Alb(X)\rar Alb(C)$. This is a contradiction since by assumption $b_1(X)=b_1(C)+2$.\\

\nind
Suppose that $F_0$ is an elliptic curve with an ordinary cusp. Then $F_0$ has a
non-constant morphism to the elliptic curve $E'$ which is the kernel of
$Alb(X)\rar Alb(C)$. Since any non-constant morphism between smooth elliptic
curves is {\'e}tale we get a contradiction by considering the induced morphism
$\ol{F_0}\rar E'$, where $\ol{F_0}$ is a smooth model of $F_0$.\\

This completes the proof of Lemma 1.\\

\nind
By Lemma 1 there is some elliptic curve $E_s$ contained in a fiber of $f$ which also maps onto
$E_0$. As defined earlier, $g:X\rar E$ is the Stein factorization of the map $X\rar E_0$. Since
the map $E_s\rar E$ is non-constant, $E$ is a smooth elliptic curve. We will
prove that this is the elliptic curve in the statement of the theorem.\\

\nind
We have a morphism $h:=f\times g:X\rar C\times E$. Then $h$ maps fibers of $f$
onto fibers of the first projection of $C\times E$ which are all isomorphic to $E$.\\

Note that the morphism $h$ is defined even when $f$ is a $C^{\infty}$-fiber bundle, since by assumption $b_1(X)=b_1(C)+2$.\\

Any rational component of any fiber of $f$ maps to a point in $C$ and $E$ both, hence maps to a point in $C\times E$.\\
Let $Z$ be the normalization of $C\times E$ in the function field of $X$. It follows that $Z$ is a projective surface. We get an induced finite morphism $h_0:Z\rar C\times E$.\\
It is clear that for any singular fiber of $f$ at least one elliptic component of this singular fiber (which exists by Lemma 1) maps onto $E$. The next result will imply that any elliptic component in a reducible fiber of $f$ maps onto $E$, so that $Z$ is obtained from $X$ by contracting rational components in fibers of $f$ to points.\\

\nind
{\bf Lemma 2.} {\it Suppose that a singular fiber $F_0$ of $f$ contains two elliptic curves $E_i$. Then the induced maps $E_i\rar E$ are both non-constant.}\\

\nind
{\bf Proof.} In the lists in \cite{NU},~\cite{O}, only the singular fiber of
type (13) in \cite{NU} contains two (smooth) elliptic curves joined by a linear
chain of $(-2)$-curves. By the proof of
Lemma 1 at least one of these two curves maps onto E, so that its image under
$h$ is not a point. Suppose $E_2$ maps to a point in $E$. If $D\subset X$ is an irreducible curve such that $f(D)=C$ then $h(D)$ is a curve. Hence the only irreducible curves in $X$ which map to points under $h$ are some components of fibers of $f$.\\
Let $Z$ be the surface defined above. Let $f_0:Z \rar C$ be the induced genus-2
fibration on $Z$. Every scheme-theoretic fiber of $f_0$ has arithmetic genus $2$
since $f_0$
is a flat map. But the image of $E_1$ in $Z$ will be a reduced curve with a unibranch
singularity. This contradicts the proof of Lemma 1.\\
This contradiction proves the result.\\

\nind
From Lemma 1 and 2 it follows that the morphism $X\rar E$ factors through $Z$.\\

\nind
For a general fiber $F$ of $f$ consider the map $pr_2\circ h:F\rar E$. By Riemann-Hurwitz formula either there are two points in $F$ which are ramified for this map with ramification index $2$ each, or a single point which is ramified with ramification index $3$. It follows that the union of these ramified points, say $R$, consists of either two irreducible curves $R_1,R_2$, or a single irreducible curve $R$. Each $R_i$ (resp. $R$) maps onto $C$ under $f$.\\ 

\nind
{\it Case 1.} The ramification locus $R$ is a cross section of $f$.\\
In this case there is a unique point in a general fiber $F$ which is ramified for the map $F\rar E$. The image of $R$ in $C\times E$, say $B$, is then a cross-section for $p_1:C\times E \rightarrow C$. Since $B$ meets every fiber of $p_1$ 
transversally the surface $Z$ is a $C^{\infty}$-fiber bundle over $C$. In this case, $X=Z$ and $f$
has no singular fibers. We will consider this case later after Step 5 and show that it cannot occur.\\

\nind
{\it Case 2.} The ramification divisor has two irreducible components $R_1,R_2$.\\
Now both $R_i$ are cross-sections for $f$ and their images $B_i$ are cross-sections for $p_1$.\\

\nind
Note that $E$ acts on $C\times E$ by translation on the second factor. Any cross-section for $p_1$ can be assumed to give the identity element in the fibers of $p_1$, considered as an abelian variety. So we will assume that $B_1$ is a standard cross-section for $p_1$.\\
In particular, $B_i^2=0$ for $i=1,2$. If $B_1\cap B_2=\phi$ then $Z$ is a $C^{\infty}$-fiber bundle.\\

\nind
{\it Case 3.} The ramification divisor is irreducible and is not a cross-section.\\

{\bf Proof of Case 2.}\\

We will first consider Case 2 assuming that there is an elliptic curve in some singular fiber which occurs with multiplicity $1$ in that fiber and use the results obtained to deal with Case 3 and the remaining (two) cases of singular fibers which contain a non-reduced elliptic curve.\\

\nind
Since the proof is somewhat long we will split it in several steps.\\

\nind
{\bf Step 1. Elimination of certain singular fibers.}\\

\nind
This step is also applicable in Case 2 and Case 3 when some elliptic curve
occurs in a singular fiber of $f$ with multiplicity 1.\\

Assume that there is a singular fiber $F_0$ of $f$ which is a union of a reduced elliptic curve $E_0$ having at most unibranch singular points and certain number of connected curves $T_1,T_2,...$ which are disjoint from each other and each of which is a union of rational curves, each $T_i$ meeting $E_0$ in only one point.\\
Each $T_i$ contracts to a normal singular point on a normal compact complex surface $X_0$ with an induced genus-$2$ fibration $f_0:X_0\rar C$. It follows that the image of $E_0$ in $X_0$ is a reduced fiber $E_0'$ having only unibranch singular points and with arithmetic genus $2$. But then $E_0'$ maps onto a smooth elliptic fiber of the map $C\times E\rar C$. This is a contradiction, as seen in the proof of Lemma 1.\\

\nind
Now we can assume in the subsequent arguments that $f$ has no such singular
fiber. This eliminates many types of singular fibers containing an elliptic
curve occuring with multiplicity $1$. From the list in \cite{NU} it now follows
that any connected component of the union of rational curves in any singular
fiber is a tree of $(-2)$-curves which can be contracted to a rational double
point. This observation will be useful later.\\

\nind
{\bf Lemma 3.} {\it The natural homomorphism $\pi_1(X)\rar\pi_1(Z)$ is an isomorphism.}\\

\nind
{\bf Proof.} This follows from the discussion in the preliminaries.\\

\nind
Lemma 3 will be useful later.\\

\nind
{\bf Step 2. The case when $R=R_1\cup R_2$ and every elliptic curve in a singular fiber occurs with multiplicity $1$.}\\

\nind
The proof in this step contains all the crucial ideas.\\ 

\nind
In view of Step 1, in the present step whenever a singular fiber contains an elliptic curve one of the following two things can happen.\\
(a) Either there are two elliptic curves in that fiber, or\\
(b) The contraction of the rational curves in that fiber does not give a (reduced) elliptic curve with at most unibranch singularities.\\

\nind
Now consider the finite morphism $h_0:Z\rar C\times E$. We denote the images of
$R_i$ in $Z$ by $R_i'$. It is clear that the ramification index of $R_i'$ over
$B_i$ is $2$ and the map $R_i'\rar B_i$ is an isomorphism since both are
cross-sections. Let $q\in B_1\cap B_2$. There is a unique point $p\in R_1'\cap
R_2'$ which maps to $q$.\\

\nind
The next result is rather technical. This result is needed for some later arguments.\\  

\nind
{\bf Lemma 4.} {\it Let $(Z,p)$ be a germ of a rational double point in
dimension~ 2. Let $\pi : (Z,p) \rar (\C^2,0)$ be a finite analytic map. Assume
that $\pi$ is ramified precisely over $\{X=0\}$ and $\{X-~Y^l=~0\}$ for some
$l\geq 1$. Further assume that there is a unique irreducible component of the
inverse image of $\{X=0\}$ in $Z$ which is ramified and has ramification index 2, and
similarly there is a unique irreducible component of the inverse image of $\{X-Y^l=0\}$ in $Z$ which is
ramified and has ramification index $2$. Then $\pi$ is Galois.}\\

\nind
{\bf Proof.} We denote the ring $\C[[X,Y]]$ by $T$.\\

\nind
If $l=1$ then the branch locus $B$ is a normal crossing divisor in the germ $(\C^2,0)$. Using the fact that $\pi_1(\C^2-B)\cong\Z\times\Z$, it is easy to see that $\pi$ is in fact an abelian extension.\\

\nind
Now we will assume that $l>1$ for the rest of the argument.\\ 

\nind
Consider the ring $T[X^{1/2}, (X-~Y^l)^{1/2}]$. This is an integral extension
of $T$ and the extension of the quotient fields is Galois with Galois group Klein's four-group $V_4$. Let $S$ be
the local ring of $Z$ at point $p$. So that we have an inclusion $T \rightarrow
S$. We wish to show that this inclusion is Galois. We have the commutative
diagram:

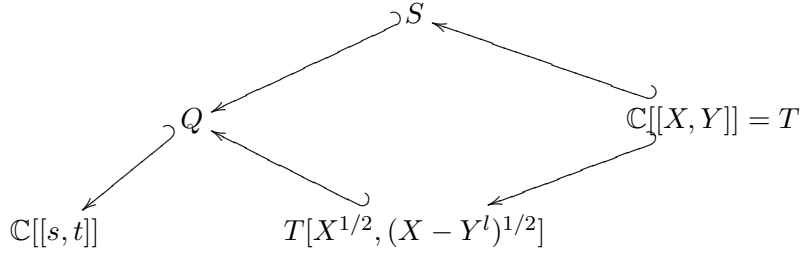
\[
\xymatrix{ & & S \ar@{_{(}->}[dl]   & & \\
& Q \ar@{_{(}->}[dl]  & & \C[[X,Y]]=T \ar@{_{(}->}[ul] \ar@{_{(}->}[dl]  & \\
\C[[s,t]] & & T[ X^{1/2}, (X-Y^l)^{1/2}] \ar@{_{(}->}[ul]  }
\]

\nind
Here $Q$ is the compositum of $S$ and $~T[ X^{1/2}, (X-~Y^l)^{1/2}]$. Since the
ramification indices over the curves $\{X=0\}$ and $\{X-Y^l=0\}$ in both the
extensions $S$ and $T[ X^{1/2}, (X-~Y^l)^{1/2}]$ are same, by Abhyankar's
lemma \cite{A} we see that
 $T[ X^{1/2}, (X-Y^l)^{1/2}] \subset Q$ is an unramified extension. Let
$u:=X^{1/2},v:=(X-Y^l)^{1/2}$. Then $T[ X^{1/2}, (X-~Y^l)^{1/2}]$ defines the
rational double point $\{u^2-v^2-Y^l=0\}$. It follows that $Q$, being
divisorially unramified over this singularity, is also a rational double point.
Now $Q$ is a ring of invariants of $\C[[s,t]]$ by a finite subgroup
of $SL_2(\C)$. Thus we have an inclusion $T \rightarrow \C[[s,t]]$ which is
{\it uniformly ramified in co-dimension one} in the sense of Griffith \cite{Gr}. This means that for any height $1$ prime ideal $P\subset T$ the ramification index of any prime ideal in $\C[[s,t]]$ lying over $P$ depends only on $P$. From
the result of (\cite{Gr}, Theorem 1.6)  we can conclude that it is Galois, with
the Galois group denoted by $G$. It is easy to see that the Galois group of the
extension $T[ X^{1/2}, (X-Y^l)^{1/2}] \subset \C[[s,t]]$ is $ \Z /l\Z$ and the
degree of the extension $T=\C[[X,Y]] \subset C[[s,t]]$ is $4l$.\\

\nind
Suppose $T\subset S$ is not Galois. Then there is a non-normal subgroup in $G$, say $H$, which corresponds to the extension  $S\subset \C[[s,t]]$ . Because $ T \rightarrow T[ X^{1/2}, (X-Y^l)^{1/2}]$ is Galois, $\Z /l\Z$ is normal in $G$ and we have the following exact sequence:

$$(1) \rightarrow \Z /l\Z \rightarrow G \rightarrow V_4 \rightarrow (1).$$

\nind
{\it Case 1.} $H$ does not intersect $ \Z /l\Z$.\\

\nind
In this case $H$ is isomorphic to $ \Z /2\Z$ or $V_4$. If $H$ is isomorphic to $
\Z /2\Z$ then the inclusion $S \subset \C[[s,t]]$ has degree 2 and it is
unramified because $S$ is not regular. This makes the map $ \C[[X,Y]] \subset S$ Galois due to the above
mentioned result of Griffith in \cite{Gr}.\\

\nind
Assume now that $H$ is isomorphic to $V_4$. We claim that $H$ is generated by
pseudo-reflections. This is seen as follows. Since $H$ is abelian its action can
be assumed to be diagonal. If $H$ is not generated by pseudo-reflections then $\sigma(s,t)=(-s,-t),~\tau(s,t)=(s,-t)$ generate
$H$. Then the ring of invariants is $\{s^2,t^2\}$. This shows that $Z$ is smooth,
a contradiction to the hypothesis.\\

\nind
{\it Case 2.} Thus we can assume that $H \cap \Z /l\Z$ is $\Z /m\Z$ for some $m \leq l$.\\

\nind
We denote the quotient of $ \C[[s,t]]$ by $\Z /m\Z $ by $W$. Then $W$ is defined
by the equation $Z_1Z_2=Z_3^m$ where $Z_1=s^m$, $Z_2=t^m$ and $Z_3=st$ are the
invariants of the action of $\Z /m\Z $ on $ \C[[s,t]]$ . As elements of
$\C[[s,t]]$ we have $u= (s^l+ t^l)/2$,  $v=(s^l-t^l)/2$ and $Y=st$. Because
$Y\in T$ the element $st$ is invariant under the action of $H$ on $W$. Again,
the following two short exact sequences are possible

$$(1) \rightarrow \Z /m\Z \rightarrow H \rightarrow \Z /2\Z \rightarrow (1)$$

$$(1) \rightarrow \Z /m\Z \rightarrow H \rightarrow V_4 \rightarrow (1).$$

\nind
In the first case let $\sigma$ be the generator of $\Z /2\Z$. Then $
\sigma(Z_3)=Z_3$ and because $Z_1Z_2=Z_3^m$, $\sigma(Z_1Z_2)= Z_1Z_2$ the
element $\sigma$ acts on $Z_1$ and $Z_2$ in the following two possible ways:
$\sigma(Z_1)=UZ_1$ and $\sigma(Z_2)=Z_2/U$ or $\sigma(Z_1)=UZ_2$ and
$\sigma(Z_2)=Z_1/U$, where $U$ is a unit. If $\sigma(Z_1)=UZ_1$ and
$\sigma(Z_2)=Z_2/U$ then because the order of $\sigma$ is 2, it is easy to see
that $\sigma(\sqrt{U})= 1/\sqrt{U}$ or  $\sigma(\sqrt{U})= -1/\sqrt{U}$. If
$\sigma(\sqrt{U})= 1/\sqrt{U}$ then let ${Z'}_1= \sqrt{U}Z_1$ and ${Z'}_2= Z_2/
\sqrt{U}$. Then $\sigma({Z'}_1)= {Z'}_1$ and $\sigma({Z'}_2)={Z'}_2$, showing
that $\sigma$ acts as the identity element, which is not possible as $\sigma $
generates $\Z /2\Z$. Thus, $\sigma(\sqrt{U})= -1/\sqrt{U}$. In which case
$\sigma({Z'}_1)= -{Z'}_1$ and $\sigma({Z'}_2)=-{Z'}_2$, and the invariants under
this action are ${{Z'}_1}^2$, ${{Z'}_2}^2$ and $Z_3$. This implies that the
singularity of $Z$ is an $A_{2m-1}$ singularity (since the degree of the extension $S\subset \C[[s,t]]$ is now $2m$), and hence the inclusion $S
\rightarrow \C[[s,t]]$ is unramified. Then the inclusion $\C[[X,Y]] \subset S$
is uniformly ramified, and hence Galois by the result of Griffith used earlier.
This means that $H$ is normal.\\

\nind
On the other hand, if $\sigma(Z_1)=UZ_2$ and $\sigma(Z_2)=Z_1/U$ then a similar analysis will show that $Z$ is smooth, which  is a contradiction to the hypothesis.\\

\nind
We now turn our attention to the second exact sequence. The extension $ W \rightarrow Z$ is Galois with Galois group $V_4$. Let $ \sigma$ and $\sigma'$ be the generators of $V_4$. All the elements of $V_4$ leave $Z_3=st$ fixed. Then by the same argument as in the case of the first exact sequence we can show that $W/ \sigma$ is either smooth or an $A_{2m-1}$ singularity. Further denoting $W/ \sigma$ by $W'$ we can argue that if $W'$ is an $A_{2m-1}$ singularity then $W' / \sigma'$ is either smooth or an $A_{4m-1}$ singularity. The former is a contradiction to hypothesis and the latter implies that $H$ is normal. If $W'$ is smooth then $W'/ \sigma' =Z$ is smooth because one of the invariants of the local ring of $W'$ under the action of $\sigma'$ is $Z_3=st$. Hence $ \sigma'$ acts as a pseudo-reflection on $W'$, implying that $Z$ is smooth. Thus $H$ is normal in $G$.\\

\nind
This completes the proof of Lemma 4.\\

\nind
Now with $Z$ as above, let $f_0:Z\rar C,~h_0:Z\rar C\times E$ be as before. Then we have the following result.\\

\nind 
{\bf Lemma 5.} {\it No component of any fiber of $f_0$ is ramified for the map
$h_0$.}\\

\nind
{\bf Proof.} Clearly no regular fiber of $f_0$ is ramified for $h_0$. Any component of a singular fiber of $f_0$ is an elliptic curve and it occurs with multiplicity $1$ in the corresponding scheme-theoretic fiber. Using $pr_1\circ h_0=1_C\circ f_0$ we deduce that no component of any fiber of $f_0$ is ramified for the map $h_0$.\\

\nind
{\bf Lemma 6.} {\it The inverse image of $B_1$ in $Z$ has $R_1'$ as a connected
component.}\\

\nind
{\bf Proof.} The branch locus for the map $Z\rar C\times E$ is $B_1\cup B_2$.
The inverse image of $B_1-B_2$ is clearly smooth since $B_1$ is smooth. Lemma 4
shows that for any point $p\in B_1\cap B_2$ the map $Z\rar C\times E$ is locally
Galois. By Lemma 5, it is {\'e}tale at points of $Z$ not lying in $R_1'\cup
R_2'$. This proves the assertion.\\

{\it Remark.} We need lemma 4 only for the proof of Lemma 6. Without Lemma 4 it is not clear, a priori, whether the inverse image of $B_i$ can contain an irreducible component which is not ramified for $h$. If a simpler argument can be given for Lemma 6 then the technical Lemma 4 can be avoided.\\ 

\nind
{\bf Lemma 7.} {\it Suppose that $f$ is not a $C^{\infty}$-fiber bundle. Then the degree of the map $Z \rar C \times E$ is 2. In particular, there is an action of $\Z/(2)$ on $Z$ with quotient  $C\times E$.}\\

\nind
{\bf Proof.} If $B_1,B_2$ do not intersect then again $Z$ is a
$C^{\infty}$-bundle and $f$ has no singular fibers. Hence $B_1.B_2>0$. We have
seen that $B_i^2=0$ for $i=1,2$. Now $(B_1 +  B_2)^2>0$ and $B_1+B_2$ is a nef and
big divisor. The divisor $R'_1 \cup
R'_2$ is connected because $R'_i\rar B_i$ is an isomorphism and $R'_i$ is a connected component of the inverse image of $B_i$. By Lemma 6, there is a tubular neighbourhood $T_1$ of
$R'_1$ such that its image under the map $Z \rar E \times C$ is a tubular
neighbourhood of $B_1$ and the full inverse image of $B_1$ in $T_1$ is $R'_1$, the ramification index of $R'_1$ over $B_1$ being 2. Thus the map $Z
\rar E \times C$ restricted to $T_1$ is of degree 2. Similarly there is a
tubular neighbourhood $T_2$ around $R'_2$ such that the map restricted to $T_2$
has degree 2. Thus the map $ Z \rar E \times C$ restricted to $T_1 \cup T_2$ has
degree 2. It is easy to see that the complement of $B_1\cup B_2$ in $C\times E$
is affine. Thus by a suitable application of Hartogs' result, the degree of the
map $ Z\rar E \times C$ is 2.\\
Since the map $h_0$ is finite the second assertion follows.\\

\nind
{\it Remark.} Since $Z\rar C\times E$ has degree $2$ and $R_i'$ is ramified for this map, we see that $R_i'$ is in fact the full inverse image of $B_i$.\\

\nind
Since $X\rar Z$ is a minimal resolution of singularities, this $\Z/(2)$-action extends naturally to $X$. Since the fixed points of a finite group acting on a smooth variety is a closed smooth subvariety $R_1,R_2$ are both smooth and disjoint. The other irreducible curves in $X$ fixed by this action are fiber components of $f$, for otherwise any horizontal fixed curve will be a ramified curve for the map $X\rar C\times E$.\\
This observation will be useful later when we consider the case when $R$ is a $2$-section.\\

\nind
{\bf Lemma 8.} {\it There is a fibration $ Z \rightarrow E$ with connected fibers such that there is a short exact sequence $$ (1)
\rightarrow \pi_1(R_1') \rightarrow \pi_1(Z) \rightarrow \pi_1(E) \rightarrow
(1)$$ where $E$ is the elliptic curve in the above discussion.}\\

\nind
{\bf Proof.} We choose the obvious map $p_2 \circ h_0 : Z \rightarrow
E$. We have seen above that $R_1'$ is the full inverse image of $B_1$, so that $R_1'$ is set-theoretically the full fiber of $p_2\circ h_0$.\\
We claim that the fibers of $p_2 \circ h_0$ are connected. To see this we note that the degree of $h_0$ is $2$. If a general fiber of this map is not connected then we consider the Stein factorization $\varphi:Z\rar E'$ of the map $Z\rar E$. Let $E_1$ be an elliptic curve in some singular fiber of $f$. Then the map $\varphi: E_1\rar E'$ is dominant. Hence $E'$ is also an elliptic curve. Consequently, the map $E'\rar E$ is {\'e}tale. There are at least two points in $E'$ over every point in $E$. It follows that the inverse image of $B_1$ contains at least two curves. This contradiction proves the claim.\\
Using the observations that $R_1'$ has ramification index $2$ over $B_1$ and $B_1$ is a fiber of the map $C\times E$ we infer that $R_1'$ occurs with multiplicity $2$ in the corresponding fiber. Consider the induced fibration $X\rar E$. The corresponding fiber of this fibration is $G_0:=R_1 \cup
\Delta$, where $\Delta$ is a disjoint union of simply-connected curves.\\

\nind
{\em Claim.} $G_0$ has multiplicity $1$, i.e. the gcd of the multiplicities of
the irreducible components of $G_0$ is $1$.\\

We first prove that under the hypothesis of Step 2 every singular fiber has a rational 
curve as an irreducible component. Indeed, if there is no rational curve in a singular fiber
then the singular fiber will be an elliptic curve with a node or a union of two
smooth elliptic curves intersecting transversally at one point. The
ramification divisor $R$ can meet such fiber only at the point
of singularity of the fiber because a non-constant map between elliptic curves
is unramified. Thus both $R_1$ and $R_2$ will have to pass through the point of
singularity of the fiber. But the intersection of $R$ with a fiber is atmost 2.
This is a contradiction.\\

Now assume that the claim is not true, then the multiplicity of $G_0$ is $2$
since $R_1$ occurs with multiplicity $2$. Let $q\in E$ be the image of $G_0$.
Let $U$ be a small disc in $E$ with center $q$ and let $t$ be a local parameter
at $q$ in $U$. Let $N$ be the inverse image of $U$ in $X$. We consider another
disc $\wt U$ with local parameter $\tau$ and the holomorphic map $\wt U\rar U$
sending $\tau\rar \tau^2$. Let $\wt N$ be the normalization of the fiber product
$N\times_U \wt U$. Then $\wt N\rar N$ is proper and {\'e}tale of degree $2$. The
inverse image of $\Delta$ splits into a disjoint union of simply-connected
curves. Let $\wt R_1$ be the inverse image of $R_1$ in $\wt N$. We can contract
the simply-connected curves to normal singular points on a normal surface $\wt
N_0$ which has induced {\'e}tale map $\wt N_0\rar Z$. The image of $\wt R_1$ in
$\wt N_0$ is therefore smooth. It is a full fiber of the map $\wt N_0\rar \wt U$
and occurs with multiplicity $1$ by construction. This easily implies that the
local ring of $\wt N_0$ at any of the singular points is regular. This
contradiction proves the claim.\\

\nind
Now by (\cite{X}, Lemma 2) we deduce the following exact sequence:
$$\pi_1(R_1') \rightarrow \pi_1(Z) \rightarrow \pi_1(E)\rightarrow (1).$$ 

\nind
The map $\pi_1(R_1') \rightarrow \pi_1(Z)$ is injective because
the following diagram commutes

$$\begin{CD}
\pi_1(R_1') @>>> \pi_1(Z)\\
@VlVV            @VVV                 \\
\pi_1(C)    @>m>>    \pi_1(E \times C)
\end{CD}
$$
where $l$ is an isomorphism and $m$ is injective. This proves the lemma.\\

\nind
{\bf Lemma 9.} {\it The fundamental group of $X$ is $\pi_1(C) \times
\pi_1(E)$}\\

\nind
{\bf Proof.} We have the following diagram

$$\begin{CD}
(1) @>>> \pi_1(R_1')@>>> \pi_1(Z) @>>> \pi_1(E) @>>> (1)\\
@.          @.            @Ab_*AA                   \\
@. \pi_1(F) @>>> \pi_1(X) @>>> \pi_1(C).@>>> (1)
\end{CD}
$$

\nind
Here $b_*$ is an isomorphism. The lower row in the diagram splits, using a
splitting $\sigma:\pi_1(C)\rar\pi_1(X)$, because $R_1\rar C$ is a cross-section
of $f$. We denote the image of $\pi_1(F)$ inside $\pi_1(X)$ by $I$. The image of
$ \sigma$ coincides with the image of $\pi_1(R_1')$ under $b_*^{-1}$ making it a
normal subgroup of $\pi_1(X)$. Since the intersection of $I$ and this normal
subgroup is trivial we have $\pi_1(X)= I \cdot \pi_1(R_1)$. The upper row shows
that $\pi_1(Z)/\pi_1(R_1')=\pi_1(E)$ thus $I=\pi_1(E)$. Since $I,\pi_1(C)$ are
both normal and have a trivial intersection they commute with each other.\\

\nind
This completes the description of $\pi_1(X)$ when $R=R_1\cup R_2$ and any elliptic curve in a singular fiber occurs with multiplicity $1$.\\

\nind
{\bf Step 3. The case when $R=R_1\cup R_2$ and some singular fiber $F_0$ contains an elliptic curve with multiplicity $2$.}\\

\nind
From \cite{NU} there are only two types of singular fibers which contain
non-reduced elliptic curves, viz. those of numerical type (12) on pages 155 and
159 of \cite{NU}. In each of these cases the fiber is a union of a smooth
elliptic curve (occurring with multiplicity $2$) and a disjoint union of at most
two trees of smooth rational curves which contract to rational double points in
$Z$.\\

We will briefly indicate the proof by using arguments from Step 2.\\

\nind
We have the exact sequence
$$\pi_1(F)\rar\pi_1(X)\rar\pi_1(C)\rar (1).$$ 
We can assume that $F$ lies in a tubular neighborhood $T$ of $F_0$. Clearly, $\pi_1(T)\cong\pi_1(F_0)\cong\Z\oplus\Z$. This gives the exact sequence 
 $$\Z\oplus\Z\rar\pi_1(X)\rar\pi_1(C)\rar (1).$$
Since $b_1(X)=b_1(C)+2$ we see that the first map in this sequence is an injection and the image is a normal subgroup of $\pi_1(X)$. The sequence splits, as before, using one of the cross-section, say $R_1$. Hence $\pi_1(X)\cong \pi_1(C) \cdot(\Z\oplus\Z)$.\\
Next, we again consider the morphism $X\rar C\times E$, inducing a morphism $Z\rar C\times E$. Under this map the factor $\Z\oplus\Z$ maps injectively to $\pi_1(E)$ and the map on the second factor $\pi_1(C)$ is an isomorphism. Using normality of $\Z\oplus\Z$ in $\pi_1(X)$ we deduce that $\pi_1(X)\cong \pi_1(C) \times(\Z\oplus\Z)$ as follows. Let $\alpha\in \Z\oplus\Z,\beta\in\pi_1(C)$. The images of these elements in $\pi_1(E),\pi_1(C)$ in $\pi_1(E\times C)$ commute. Hence $\alpha,\beta^{-1}\alpha\beta$ map to the same element in $\pi_1(E\times C)$. But the map $(\Z\oplus\Z)\rar\pi_1(E)$ is injective. This shows that $\alpha,\beta$ commute.\\

\nind
This completes the proof of the part (1) of the theorem when $R$ is a union of two cross-sections.\\ 

{\bf Proof of Case 3.}\\

\nind
{\bf Step 4. The case when $R$ is irreducible, not a cross-section, and there is at least one singular fiber containing a reduced elliptic curve.}\\

\nind
Again, let $f: X \rightarrow C$ be a genus-2 fibration satisfying the properties
mentioned above with $R$ irreducible and not a cross-section. In this
sub-section we prove that the degree of the map $X \rightarrow E \times C$ is 2
and discuss the structure $\pi_1(X)$. These observations will help us verify the
Shafarevich conjecture for $X$.\\

{\bf Lemma 10.} {\it The map $X \rightarrow E \times C$ has degree 2.}\\

{\bf Proof.} We have a map $f|_R:R \rightarrow C$ of degree $2$. We consider the desingularized fiber
product $Y=R \times_C X$. The fibration $f': Y \rightarrow R$ is a genus-2
fibration with connected fibers. Since $f$ has a fiber, say $F_s$, containing a reduced elliptic component, say $E_s$, we can see that the same thing is true for $f'$.\\
We briefly indicate a proof of this. If the map $R\rar C$ is ramified at a point
in $R$ where $R$ meets $F_s$ then an easy argument shows that $E_s$ is ramified
for the map $Y\rar X$. Since the degree of the map $Y\rar X$ is $2$ the inverse
image of $E_s$ in $Y$ contains an elliptic curve. On the other hand, if the map
$R\rar C$ is {\'e}tale at points in $R\cap F_s$ then $Y\rar X$ is {\'e}tale over
$F_s$. Again, the inverse image of $E_s$ in $Y$ contains an elliptic curve.\\ 
Hence $f'$ cannot be a $C^{\infty}$-fiber bundle. Since the image of $H_1(F,{\bf R})\rar H_1(X,{\bf R})$ has image of rank $2$, the same thing is true for $Y$. Hence $b_1(Y)=b_1(R)+2$.\\

There is a map $l:Y \rightarrow E$ given by the
composition $Y \rightarrow X \rightarrow E$. Now $f'$ has a tautological cross-section $R_1'$ such that there are two cross-sections,$R_1',R_2'$, of $f'$ which lie over $R$. By the previous arguments the map $Y \rightarrow R \times E$ is
a degree 2 map. Now we consider the commutative diagram:
$$\begin{CD}
Y  @>>> X\\
@VVV            @VVV                 \\
R \times E    @>m>>    C \times E
\end{CD}
$$

By the proof of Step 2, the degrees of the maps $Y \rightarrow R \times E$, $R \times E \rightarrow C
\times E$ and $Y \rightarrow X$ are all 2. Thus the map $X \rightarrow C \times E$
has degree 2.\\

We will need the following general result.\\

{\bf Lemma 11.} {\it Let $f:W\rightarrow V$ be a finite surjective analytic map
between irreducible normal complex spaces. Suppose that for some point $p\in V$
there is a unique point $q\in W$ lying over $p$. Then the induced natural
homomorphism
$\pi_1(W)\rightarrow\pi_1(V)$ is a surjection.}\\

{\bf Proof.} Let $B\subset V$ be the branch locus. Then
$f:W-f^{-1}(B)\rightarrow V-B$ is a topological covering of finite degree, say
$d$. Hence the homomorphism $\pi_1(W-f^{-1}(B))\rightarrow\pi_1 (V-B)$ is
injective and the image has index $d$. Both the homomorphisms
$\pi_1(W-f^{-1}(B))\rightarrow\pi_1(W),~\pi_1(V-B)\rightarrow \pi_1(V)$ are
surjections because $V,W$ are normal. It follows that the image of the homomorphism
$\pi_1(W)\rightarrow\pi_1(V)$ has index at most $d$. Let $H$ be the image of
this homomorphism. There is a  connected covering $\tilde V\rightarrow V$ such
that $\pi_1(\tilde V)=H$. If the above index is $>1$ then there are at least two distinct points in $\tilde
V$ lying over $p$. By covering space theory we have a lift $W\rightarrow\tilde
V$, which will also be a finite surjective analytic map. But $W$ has a unique
point $q$ lying over $p$. This is a contradiction.\\
Hence $\pi_1(W)\rightarrow\pi_1(V)$ is a surjection.\\

By the previous step, $\pi_1(Y)\cong \pi_1(R)\times\pi_1(E_0)$, where $E_0$ is an elliptic curve. If the map $f|_R$ is unramified then the map $\pi_1(Y) \rightarrow \pi_1(X)$ is
injective. Thus we have $\pi_1(R) \times \pi_1(E_0)$ as an index 2 subgroup of
$\pi_1(X)$. If the map $f|_R$ is ramified then we have the following Lemma.\\

\nind
{\bf Lemma 12.} {\it Assume that $f|_R$ is ramified. Then $\pi_1(X)= \Gamma_1
\times \pi_1(E_0)$ where $\Gamma_1$ is a surjective image of $\pi_1(R)$.}\\
{\bf Proof.} Consider the commutative diagram:
$$
\xymatrix{
\pi_1(Y) \ar[d]^b \ar[r]^a  &\pi_1(X) \ar[d]^c\\
\pi_1(R \times E) \ar[r]^d  &\pi_1(C\times E)}
$$

The maps $a$ and $d$ are surjective by Lemma 11. The map $R \times E
\rightarrow C \times E$ is given by $f|_R \times id$. The morphism $Y\rar R$ may
still have some singular fiber containing a non-reduced elliptic curve. Hence we
have an isomorphism $\pi_1(Y)\cong \pi_1(R)\times \pi_1(E_0)$. Under the map
$b$ the factor $\pi_1(E_0)$ maps injectively to $\pi_1(E)$. Let $\Gamma_1$ and
$\Gamma_2$ be the images
of $\pi_1(R)$ and $\pi_1(E_0)$ in $\pi_1(X)$ under the map $\pi_1(Y)\rar\pi_1(X)$. We claim that
$\Gamma_1 \cap \Gamma_2 = {1}$. This follows from the observation that the images of $\Gamma_1,\Gamma_2$ are mapped respectively to $\pi_1(C),\pi_1(E)$ by $c$ in the group $\pi_1(C\times E)$.\\

\nind
{\bf Step 5. The case when $R$ is irreducible and every singular fiber contains a non-reduced elliptic curve.}\\

We again consider the desingularized fiber product $Y: X\times_C R$. First of all, $Y$ may not be relatively minimal. Since the fundamental group does not change by contraction of a $(-1)$-curve, for the sake of simplicity of exposition, we will assume that $f':Y\rar R$ is relatively minimal.\\

{\it Claim.} $f'$ is not a $C^{\infty}$-bundle.\\

To see this we will look closely at the two types of singular fibers $F_s$ which contain an elliptic curve $E_s$ occuring with multiplicity $2$.\\
Suppose $R$ meets $E_s$. Since the degree of the map $R\rar C$ is $2$, the intersection $E_s\cap R$ is a single point and the intersection is transverse. But then the map $h|_{E_s}:E_s\rar E$ is ramified. This is not possible. Hence $R$ does not meet $E_s$.\\
Suppose that $R\rar C$ is {\'e}tale at all the points in $F_s\cap R$. Then clearly the inverse image of $E_s$ in $Y$ is a union of elliptic curves. Hence $f'$ is not a bundle.\\
Now assume that $R\rar C$ is ramified at some point, say $p$, in $F_s\cap R$. From the nature of $F_s$ we see that $p$ lies in an irreducible rational component $F_1$ of $F_s$. If this component occurs with multiplicity $2$ in $F_s$ then again the map $Y\rar X$ is {\'e}tale in a neighborhood of $E_s$. Hence $f'$ cannot be a bundle.\\
Assume finally that $F_1$ occurs with multiplicity $1$. In one type of singular fiber $E_s$ and $F_1$ are connected by rational components each of which occurs with multiplicity $2$ in $F_s$. Hence $Y\rar X$ is {\'e}tale in a neighborhood of $E_s$, proving that $f'$ is not a bundle.\\
Suppose $F_1$ meets $E_s$. Then their intersection is a single point with a transverse intersection. Then we see that $F_1$ is ramified for the map $Y\rar X$ and the inverse image of $E_s$ in $Y$ is a smooth irreducible curve $\wt{E_s}$, the map $\wt{E_s}\rar E_s$ has degree $2$ and there is a unique point in $\wt{E_s}$ which is ramified for the map to $E_s$ with ramification index is $2$. This contradicts the Riemann-Hurwitz formula.\\
This contradiction proves the claim.\\ 

By the claim above, $f'$ is not a $C^{\infty}$-fiber bundle. Then $b_1(Y)=b_1(R)+2$. Now the ramification divisor for the morphism $Y\rar R\times E$ splits into two cross-sections (one of which is the tautological cross-section). We can now use Step 2 or Step 3 for $Y$. Hence $\pi_1(Y)\cong\pi_1(R)\times(\Z\oplus\Z)$. If $R\rar C$ is {\'e}tale (of degree $2$) then $Y\rar X$ is {\'e}tale (of degree $2$). Hence $\pi_1(R)\times (\Z\oplus\Z)$ is a subgroup of index $2$ of $\pi_1(X)$.\\

Suppose that $R\rar C$ is ramified. Then $\pi_1(R)\rar\pi_1(C)$ is onto by Lemma 11. In this
case we deduce that $\pi_1(X)\cong \Gamma \times (\Z\oplus\Z)$, where
$\Gamma$ is the image of $\pi_1(R)\rar\pi_1(X)$.\\

{\bf Proof of Case 1.}\\

Now we consider the remaining case when $R$ is a cross-section.\\
We will show that this case cannot occur.\\

The degree of the map $R\rar C$ is $3$. In this case the branch curve $B$ is a cross-section for the map $C\times E\rar C$. Hence the map $Z\rar C$ is a $C^{\infty}$-bundle, $X=Z$ and $f$ is a $C^{\infty}$-bundle. Since $B$ is smooth, by a local analysis we see easily that its inverse image in $X$ is smooth. Since $R$ is the ramification divisor in $X$ we deduce that $R$ is the full inverse image of $B$ with ramification index $3$. Now $C$ cannot be rational, otherwise
 by Teichm{\"u}ller theory $X\cong C\times F$ which contradicts the assumption $b_1(X)=b_1(C)+2$. The long exact homotopy sequence for the bundle $f$ gives a short exact sequence
$$(1)\rar \pi_1(F)\rar\pi_1(X)\rar\pi_1(C)\rar(1).$$
Using the cross-section $R$ we have a splitting $\pi_1X)\cong\pi_1(F)\cdot\pi_1(C)$, where $\pi_1(F)$ is normal by the above short exact sequence.\\
Since $R$ is the full inverse image of $B$, the proof of Lemma 8 gives a short exact sequence 
$$(1)\rar\pi_1(R)\rar\pi_1(X)\rar\pi_1(E)\rar (1).$$ 
This shows that $\pi_1(R)\cong\pi_1(C)$ is a normal subgroup of $\pi_1(X)$. This finally shows that $\pi_1(X)\cong \pi_1(F)\times\pi_1(C)$. But then $b_1(X)=b_1(C)+4$, a contradiction.\\

This completes the proof of part (1) of the theorem.\\

{\bf Proof of part (2)}.\\

Suppose now that $X \rar C$ is a genus-2 fibration with no singular fibers.\\
If $C$ is rational then using Teichm{\"u}ller theory we know that $f$ is a trivial bundle.\\

\nind
Suppose that genus of $C$ is $1$. Again by Teichm{\"u}ller theory all the fibers of $f$ are mutually isomorphic. In this case it is known that after taking a suitable finite unramified cover $\tilde C\rar C$ the fiber product $X\times_C\tilde C$ is isomorphic to $\tilde C\times F$.\\

This completes the proof of the theorem.\\

\nind
{\it Remark.} In the (2005) Ph.D. thesis of G. Karado{\u{g}}an (\cite{K}) at the Middle East Technical University it is proved that when $f$ is a $C^{\infty}$-bundle and $b_1(X)=b_1(C)+2$, there is another morphism $\varphi:X\rar E$, where $E$ is an elliptic curve, such that $\varphi$ has exactly two multiple fibers of multiplicity $2$ each (the reduced fiber being smooth) and no other singular fibers.\\

\section{Proof of Corollary 1. The Shafarevich Conjecture}

\nind
We now prove the Corollary 1 stated in the Introduction. The following is a well-known question raised by I.R. Shafarevich (\cite{Sh}, Chapter IX, \S 3).\\

\nind
{\bf Conjecture.} {\it Let $X$ be a smooth projective variety over $\C$. Then the universal covering space $\wt{X}$ of $X$ is holomorphically convex.}\\

\nind
Recall that a complex manifold $U$ is holomorphically convex if for any sequence
of points $u_n$ in $U$ without a limit point there is a holomorphic function $f$
on $U$ such that the sequence $|f(u_n)|$ is unbounded. A compact complex
manifold is trivially holomorphically convex. Any non-compact connected Riemann surface is
Stein by a well-known result due to R. Narasimhan \cite{N}, hence
holomorphically convex. In \cite{GS} this conjecture was verified for all
smooth projective surfaces which are not of general type. Here we verify the
conjecture for genus-2 fibrations $f: X \rightarrow C$ fibrations such that the
map $b_1(X)> b_1(C)$.\\

\nind
Again by Beauville's result mentioned before and the proposition in the preliminaries, if $b_1(X)=b_1(C)+4$ then $X$ is isomorphic to the product $C\times F$. Since the universal cover of $C\times F$ is the product of the universal covers of $C,F$, which is holomorphically convex, the result follows.\\

\nind
In view of this we will now assume that $b_1(X)=b_1(C)+2$.\\

\nind
For later use we state here the following basic result due to R. Remmert and K. Stein. (\cite{RS}, Satz 8).\\

\nind
{\it Let $\tau:W\rar V$ be a proper surjective holomorphic mapping between
normal complex spaces of same dimension. Then $V$ is holomorphically convex if and only if $W$ is
holomorphically convex.}\\ 

\nind
{\bf Case 2.} {\it $R$ is a union of two cross-sections.}
We have shown in this case that either\\

(1) $\pi_1(X) \cong \pi_1(C) \times \pi_1(E)$, where $E$ is an elliptic curve and the induced map on the fundamental groups corresponding to a suitable morphism $X\rar C\times E$ is this isomorphism, or\\

(2) $\pi_1(X)\cong\pi_1(C)\times\pi_1(E_0)$, where $E_0$ is an elliptic curve (occuring in a singular fiber).\\

Consider the case (1).\\

We have therefore the commutative diagram:

$$\begin{CD}
\tilde{X}            @>h>>       U \times \C \\
     @Vv_1VV                           @Vv_2VV                 \\
       X         @>f \times g>>     C \times E
\end{CD}
$$

\nind
Where $\tilde{X}$ is the fiber-product of the universal cover $U \times \C$
of $E \times C$ and $X$. The space $\tilde{X}$ is the universal cover of $X$. Now $\tilde{X}$ is holomorphically
convex because $\tilde{X} \rar U \times \C$ is a proper map with generically
finite fibers and $U \times \C$ is easily seen to be
holomorphically convex.\\

Next consider (2).\\

In this case the homomorphism $\pi_1(X)\cong\pi_1(C)\times\pi_1(E_0)\rar\pi_1(C\times E)$ is an isomorphism on the first factor and an injection with image of finite index on the second factor.\\

Again, the universal cover $\tilde X$ of $X$ is a connected component of the pull back of $U\times \C$ and the map $\tilde X\rar U\times\C$ is proper with finite fibers. Using the Remmert-Stein theorem we see that $\tilde X$ is holomorphically convex.\\

\nind
{\bf Case 3.} {\it $R$ is irreducible and not a cross-section.}\\

The map $f:X \rightarrow C$ restricts to the map $f|_R: R \rightarrow C$. 
Hence we can consider the desingularized fiber product $Y = R \times_C X$. If $f|_R$ is
unramified then the map $Y \rightarrow X$ is {\'e}tale. This implies that the
universal cover of $Y$ is the universal cover of $X$. We know the conjecture to
be true for $Y$ and hence it is true for $X$.\\

Now suppose that $f|_R$ is ramified. Then the map $Y \rightarrow X$ is of
degree 2, hence $\pi_1(Y) \rightarrow \pi_1(X)$ is either an injection with image of index $2$, or $\pi_1(X)\cong\Gamma_1\times\pi_1(E)$ such that there is a surjection $\pi_1(R)\rar\Gamma_1$.\\

In the first case, using the holomorphic convexity of the universal cover of $Y$ we deduce that for the universal cover of $X$.\\

In the second case, pull-back of the universal cover $\tilde X$ of $X$ to $Y$ is of the form $U\times\C$, where $U\rar R$ is a suitable cover.\\
Since any Riemann surface is holomorphically convex we again deduce that $\tilde X$ is holomorphically convex.\\

\nind
{\bf Case 1.} {\it $R$ is irreducible and a cross-section.}\\

We have shown that this case cannot occur.\\

\nind
This completes the proof of Corollary 1 mentioned in the introduction.\\

\nind
In \cite{Na} the author proves this result using powerful analytic methods. We
believe that our proof, which describes the fundamental group of $X$ explicitly,
sheds more light on this result.\\

\nind
In \cite{G} the following result is proved.\\

\nind
{\it If $X$ is a smooth projective surface for which the Shafarevich Conjecture is true then $\pi_2(X)$ is a free abelian group.}\\

\nind
Hence we get the Corollary 2 in the introduction as a consequence of our
Theorem.\\

\vspace{10mm}
\nind
R.V. Gurjar, School of Mathematics, Tata Institute of Fundamental Research, Homi-Bhabha Road, Mumbai-400005, India.\\
e-mail: gurjar@math.tifr.res.in\\

\nind
Sagar Kolte, School of Mathematics, Tata Institute of Fundamental Research, Homi-Bhabha Road, Mumbai-400005, India.\\
e-mail: sagar@math.tifr.res.in

\end{document}